\documentclass[twoside]{article}
\usepackage{graphicx}
\usepackage{amsmath}
\usepackage{amssymb}
\usepackage{amscd}
\usepackage{hhline}
\newtheorem{theorem}{Theorem}

\newtheorem{theoremb}{Theorem}
\newtheorem{theoremc}{Theorem}
\newtheorem{theoremd}{Theorem}
\newtheorem{theoreme}{Theorem}
\newtheorem{dfn}[theoremb]{Definition}
\newtheorem{rk}[theoremc]{Remark}
\newtheorem{cor}[theoremc]{Corollary.}
\newtheorem{lem}[theoreme]{Lemma}
\newtheorem{examp}[theoremd]{Example}
\newtheorem{prop}[theoreme]{Proposition}
\newenvironment{proof}[1][Proof]{\textbf{#1.} }{\qed \vspace{5pt}}
\newenvironment{Proof}[1]{\textbf{#1.} }{\qed \vspace{5pt}}

\newcommand\bib[1]{\bibitem[#1]{#1}}

\newcommand\abz{\hspace{12pt}}
\newcommand\qed{\phantom{\underline{y}}\hfill\hfill$\square$}

\renewcommand\a{\alpha}
\renewcommand\b{\beta}

\newcommand\C{{\mathbb C}}

\renewcommand\d{\delta}

\newcommand\hps{\hskip-16pt . \hskip2pt}

\renewcommand\l{\lambda}

\newcommand\La{\Lambda}

\newcommand\oo{\omega}
\newcommand\op[1]{\mathop{\rm #1}\nolimits}
\newcommand\ot{\otimes}
\newcommand\p{\partial}

\newcommand\po{$\!\!\!{\text{\bf.}}$ }
\newcommand\poo{$\hspace{-8pt}{\text{\bf.}}$\ \,}

\newcommand\R{{\mathbb R}}

\newcommand\te{\theta}
\newcommand\tg{\bar g}
\newcommand\ti{\tilde}

\newcommand\vp{\varphi}
\newcommand\we{\wedge}
\newcommand\x{\xi}

\newcommand\Z{{\mathbb Z}}
\newcommand{\weg}[1]{}

\makeatletter
\renewcommand{\@oddhead}{\hfil Involutive symbolic PDEs\hfil}
\renewcommand{\@evenhead}{\hfil Boris Kruglikov, Valentin Lychagin\hfil}
\makeatother

\begin{document}

 \title{Spencer $\d$-cohomology, restrictions, \\
characteristics and involutive symbolic PDEs}
 \author{Boris Kruglikov, Valentin Lychagin}
 \date{}
\maketitle

 \vspace{-14.5pt}
 \begin{abstract}
We generalize the notion of involutivity to systems of
differential equations of different orders and show that the
classical results relating involutivity, restrictions,
characteristics and characteristicity, known for first order
systems, extend to the general context. This involves, in
particular, a new definition of strong characteristicity. The
proof exploits a spectral sequence relating Spencer
$\d$-cohomology of a symbolic system and its restriction to a
non-characteristic subspace.
 \footnote{MSC numbers: 35N10, 58A20, 58H10; 35A30. Keywords:
Spencer cohomology, symbolic system, restriction, involutivity,
characteristics.}%
 \end{abstract}

\section*{Introduction}\label{S0}

 \hspace{13.5pt}
This paper concerns some algebraic aspects of systems of
differential equations. We will investigate their systems of
symbols $g_k\subset S^kT^*\ot N$, where $T$ and $N$ are finite
dimensional vector spaces, representing the spaces of independent
and dependent variables respectively. The collection
$g=\{g_k\}_{k\ge0}$ will be called a symbolic system and we do not
require that it is generated in one particular order (more details
will be given in \S\ref{S1} below).

Spencer $\d$-cohomology groups $H^{i,j}(g)$ are algebraic
invariants of such structures, important in the study of formal
integrability of PDEs \cite{S}. Let $W\subset T$ be a subspace and
$V^*=\op{ann}(W)$. From the exact sequence
 $$
0\to V^*\hookrightarrow T^*\to W^*\to 0\eqno (\dag)
 $$
we get the restricted symbolic system $\tg_k\subset S^kW^*\ot N$.

A subspace $V^*\subset T^*$ is called non-characteristic
(\cite{S}) if no non-zero element $\oo\in g_k$ restricts to zero
on $W$: $\oo(\x_1,\dots,\x_k)=0\,\forall\x_i\in
W\Rightarrow\oo=0$, where $k=r_\text{min}(g)$ is the minimal order
of the system $g$ (\S\ref{S2}). This classical definition leads
however to the following confusion for higher-order systems
($k>1$). Consider a subspace $g_k\subset S^kT^*$ of codimension 1,
$k\ge2$, and let $g$ be the corresponding symbolic system. Then by
dimension reasons any one-dimensional subspace $V^*$ is
characteristic (see Example \ref{ex4}). Identifying
such subspaces with projectivized covectors we would conclude that
every covector is characteristic for one scalar PDE (not ODE:
$\dim T>1$) of order $k>1$, which can't be true.

Thus we risk changing the standard terminology and call
non-characteristic (by Spen\-cer and al) subspaces $V^*$ strongly
non-characteristic (with opposite being weakly characteristic; so
in the above example all covectors are weakly characteristic, but
not characteristic); more motivations for this will come below.

Let us denote
 $$
\Upsilon^{i,j}=\bigoplus_{r>0}S^rV^*\ot\d(S^{i+1-r}W^*\ot\La^{j-1}W^*)\ot
N,\quad \Theta^{i,j}=\bigoplus_{q>0}\Upsilon^{i,q}\ot\La^{j-q}V^*,
 $$
where $\d$ is the Spencer operator. Then we have:
$\Upsilon^{i,0}=\Upsilon^{0,j}=0$ and
$\Theta^{i,0}=\Theta^{0,j}=0$. The other terms are however
non-zero.

Let also $\Pi^{i,j}=\d(S^{i+1}V^*\ot N\ot \La^{j-1}V^*)$. Note
that $\Pi^{i,0}=0$ by definition and for $j>0$ we have:
$\Pi^{i,j}=S^iV^*\ot N\ot \La^jV^*\cap\op{Ker}\d$. In particular,
$\Pi^{0,j}=N\ot\La^jV^*$ for $j>0$. Finally denote
$S^{i,j}=\op{Im}(S^{i+j}V^*\to S^iV^*\ot S^jV^*)$.

Involutivity of a symbolic system $g$ is equivalent to vanishing
of certain $\d$-cohomology groups, see \S\ref{S3}.

 \begin{theorem}\po\label{thm1}
Let $V^*$ be a strongly non-characteristic subspace for a symbolic
system $g$. If $g$ is involutive, then its $W$-restriction $\tg$
is also involutive.

Moreover the Spencer cohomology of $g$ and $\tg$ are related by
the formula:
 $$
H^{i,j}(g)\simeq\, \bigoplus\limits_{q>0}\, H^{i,q}(\tg)
\ot\La^{j-q}V^*\oplus \d^{i+1}_{r_\text{min}(g)}\cdot
[\Theta^{i,j} \oplus \Pi^{i,j}]\oplus \d^i_0\d^j_0\cdot
H^{0,0}(\tg),
 $$
where $\d^t_s$ is the Kronecker symbol.

If $\tg$ is an involutive system of pure order
$k=r_\text{min}(\tg)=r_\text{max}(\tg)$, then $g$ is also an
involutive system of pure order $k$ and the above formula holds.
 \end{theorem}

 \begin{cor}
If $g$ is an involutive system and $V^*$ is strongly
non-charac\-te\-ristic subspace of $T^*$, then the
$\d$-differential induces the following exact sequences:
 \begin{multline*}
0\to
\d^{i+1}_{r_\text{min}(g)}\cdot S^{j,i}\ot N
\to S^{j-1}V^*\ot H^{i,1}(g)\to\dots\\
{}\dots\to V^*\ot H^{i,j-1}(g) \to H^{i,j}(g)\to
H^{i,j}(\tg)\oplus\d^{i+
1}_{r_\text{min}(g)}\cdot\Upsilon^{i,j}\to 0,
 \end{multline*}
for $i\ge r_\text{min}(g)-1$.
 \end{cor}

The implication ($g$ involutive)$\,\Rightarrow\,$($\tg$
involutive) from Theorem \ref{thm1} constitutes Guillemin's
theorem A obtained for the first order systems in \cite{G}. Our
proof is based on the technique, developed for other purposes in
\cite{KL$_2$}, which allows to generalize the statement to
arbitrary symbolic systems. The inverse Guillemin's theorem, i.e.
the implication, ($\tg$ involutive)$\,\Rightarrow\,$($g$
involutive) was not known.

The above Corollary for $k=0$ and system $g$ of the first order is
a theorem of Quillen-Guillemin \cite{G, Q} (see \S\ref{S4.5} for
details).

This paper is a generalization of the classical results known for
the first order systems (\cite{G,GK,BCG$^3$}). But it is not
straightforward. Indeed, Spencer's reformulation of Guillemin's
results for higher order holds true in the stable range $m\ge\mu$
(\S1.7-1.8 of \cite{S}), but one inevitably gets into the trouble
if adjusts the order $k\ge r_\text{min}(g)$ precisely to the
place, where involutivity of $g$ starts (see \S\ref{S6}). In
addition, Spencer's generalization of results due to Quillen and
Guillemin (loc.cit.) contains mistakes (see Remark in \S\ref{S2}).

More arguably, it was not noticed previously that the definition
of characteristic subspace has two meaningful generalizations to
the case of higher order: one important for the Cauchy problem
(standard one adapted to restrictions as in \cite{S}) and the
other one important in the study of characteristics, which we call
strong characteristicity. Namely, we call $V^*$ strongly
characteristic for $g$ if $\exists\,\oo\in g_k\setminus\{0\}$ such
that the directional derivative $\d_\x\oo=0\,\forall\x\in W$,
where $k=r_\text{max}(g)$ is the maximal order of the system.

 \begin{theorem}\po\label{thm2}
Let $g$ be an involutive system over $\C$. Then a subspace
$V^*\subset T^*$ is strongly characteristic iff it contains a
characteristic covector.
 \end{theorem}

This result is a generalization of Guillemin's theorem B from
\cite{G}, which concerns the pure order one systems -- the only
case, where the two introduced notions of characteristicity
coincide. Thus both Guillemin's theorems have analogs for higher
(and even various) order systems, but for this two different
notions of characteristicity should be imposed.

Some other results appear at the end of this paper. We also
provide a series of counter-examples showing importance of all our
hypotheses.

\section{\hps Symbolic systems}\label{S1}

 \abz
Consider the Spencer $\d$-complex:
 $$
0\to S^kT^*\ot N\stackrel{\d}\to S^{k-1}T^*\ot N\ot T^*\stackrel{\d}\to
\dots\stackrel{\d}\to S^{k-n}T^*\ot N\ot\La^nT^*\to0,
 $$
where $S^iT^*=0$ for $i<0$. The {\em first prolongation\/} of a
subspace $h\subset S^kT^*\ot N$ is
 $$
h^{(1)}=\{p\in S^{k+1}T^*\ot N\,|\,\d p\in h\ot T^*\}
 $$
Higher prolongations are defined inductively and satisfy
$(h^{(l)})^{(m)}=h^{(l+m)}$. An alternative definition is:
$h^{(l)}=S^lT^*\ot h\cap S^{k+l}T^*\ot N$.

 \begin{dfn}\po
Symbolic system is a sequence of subspaces $g_k\subset S^kT^*\ot
N$, $k\ge0$, with $g_0=N$ and $g_k\subset g_{k-1}^{(1)}$.
 \end{dfn}

With every such a system we associate its Spencer $\d$-complex of order
$k$:
 $$
0\to g_k\stackrel{\d}\to g_{k-1}\ot T^*\stackrel{\d}\to
g_{k-2}\ot\La^2T^*\to\dots\stackrel{\d}\to g_{k-n}\ot\La^nT^*\to0.
 $$
The cohomology group at the term $g_i\ot\La^jT^*$ is denoted by
$H^{i,j}(g)$ and is called the Spencer $\d$-cohomology of $g$.

Note that $g_k=S^kT^*\ot N$ for $0\le k<r$ and the first number
$r=r_\text{min}(g)$, where the equality is violated is called the
minimal order of the system. Actually the system has several
orders:
 $$
\op{ord}(g)=\{k\in\Z_+\,|\,g_k\ne g_{k-1}^{(1)}\}.
 $$
Note that multiplicity of an order $m(r)=\dim g_{r-1}^{(1)}/g_r$
is equal to $\dim H^{r-1,1}(g)$. Hilbert basis theorem implies
finiteness of the set of orders:
 $$
\op{codim}(g):=\dim H^{*,1}(g)=\#\op{ord}(g)<\infty.
 $$
If $r_\text{max}$ is the maximal order of the system, then
$g_{k+1}=g_k^{(1)}$ for $k\ge r_\text{max}$.

Denote by $\tg$ the image of the restriction map $g\to SW^*\ot N$.

 \begin{prop}\po
For any subspace $W\subset T$ the restriction $\tg$ is a symbolic
system.
 \end{prop}

 \begin{proof}
This follows from naturality of the $\d$-differential.
 \end{proof}

\section{\hps Characteristics}\label{S2}

 \abz
A covector $v\in {}^\C T^*\setminus\{0\}$ is called (complex)
characteristic for $g_k$ if $v^k\ot w\in g_k^\C$ for some $w\in
{}^\C N\setminus\{0\}$ (in this paper characteristics will be
considered only over the field $\C$). Clearly, if $v$ is
characteristic for $g_k$, it is characteristic for its
prolongation $g_k^{(1)}$ and vice versa. We call $v$
characteristic for a symbolic system $g$ if it is characteristic
on every level $g_k$ (or equivalently only for the level
$k=r_\text{max}(g)$). The projectivized set of all characteristic
covectors forms the characteristic variety $\op{Char}^\C(g)\subset
P^\C T^*$.

 \begin{dfn}\po
Call a subspace $V^*\subset T^*$ strongly non-characteristic for
$g_k$ if $g_k\cap V^*\cdot S^{k-1}T^*\ot N=0$. In the opposite
case, when the intersection is non-zero, let's call $V^*$ weakly
characteristic.

Call $V^*$ weakly non-characteristic if $g_k\cap S^kV^*\ot N=0$.
If the intersection is non-zero, the subspace $V^*$ will be called
strongly characteristic.
 \end{dfn}

Note that strong characteristicity implies weak characteristicity,
as well as strong non-characteristicity implies weak
non-characteristicity, but there are spaces $g_k$, for which
certain $V^*$ are (don't be confused!) simultaneously weakly
characteristic and weakly non-characteristic. Of course, then they
are neither strongly characteristic nor strongly
non-characteristic for $g_k$.

 \begin{rk}\poo
Spencer's notion of non-characteristicity (definition 1.8.1 of
\cite{S}) formally coincides with our weak non-characteristicity
(after translation from $D$-complex to the symbolic language). But
then his Theorem 1.8.1(i) (as well as preceding Theorem 1.7.3)
becomes wrong, since weak non-characteristic subspaces don't need
to be strong non-characteristic (see \S\ref{S6}). However changing
tensorial product to symmetric and shifting the index by one, the
definition turns into our strong non-characteristicity and the
subsequent statements hold. Thus importance of distinction between
characteristicities becomes apparent.
 \end{rk}

Two introduced notions of weak and strong characteristicity
coincide for first order systems, $k=1$ (the same for
non-characteristicity), but not for the case of higher orders. The
following property follows directly from the definition:

 \begin{prop}\po
A subspace $V^*$ is strongly non-characteristic for $g_k$ iff the
restriction to $W$ map $g_k\to\tg_k$ is an isomorphism. It is
weakly non-characteristic iff kernels of the maps $\d_w:g_k\to
g_{k-1}$, $w\in W$, jointly intersect only by zero.
 \end{prop}

 \begin{prop}\po\label{prop2}
If $g_k\subset S^kT^*\ot N$ is weakly/strongly non-characteristic,
then any subspace of its prolongation $g_{k+1}\subset
g_k^{(1)}\subset S^{k+1}T^*\ot N$ is such as well.
 \end{prop}

 \begin{proof}
For the weak case the statement is obvious.

Let $V^*$ be strongly non-characteristic for $g_k$. Then the
restriction $g_k\to\tg_k$ is an isomorphism. Assume that $V^*$ is
weakly characteristic for $g_{k+1}$. From the commutative diagram
 $$
  \begin{array}{ccc}
g_{k+1} & \stackrel{\d}\longrightarrow & g_k\ot T^*\\
\downarrow && \downarrow \\
\tg_{k+1} & \stackrel{\d}\longrightarrow & \tg_k\ot W^*
  \end{array}
 $$
we conclude that a non-zero element $p\in g_{k+1}$ belongs to the
kernel of the restriction map iff $\d(p)\in g_k\ot V^*$, which
implies that $V^*$ is strongly characteristic for $g_{k+1}$ and
hence strongly characteristic for $g_k$. But this yields that
$V^*$ is weakly characteristic for $g_k$, contradicting our
assumption.
 \end{proof}

Thus weak and strong non-characteristicity for a space $g_k$ are
inherited by the prolonged spaces $g_k^{(1)}$. Theorem \ref{thm2}
assures that the same is true for strong characteristicity (over
$\C$) in involutive case, but thanks to Example \ref{examp6} not
in general. As will be seen in Example \ref{examp7} the property
of $g_k$ being weakly characteristic is not hereditary upon
prolongations too.

 \begin{dfn}\po\label{def3}
A subspace $V^*\subset T^*$ is called weakly or strongly
non-characteristic for a symbolic system $g$ if this requirement
holds for $g_{r_\text{min}(g)}$ and hence for every $g_k$ with
$k\ge r_\text{min}(g)$.

Call $V^*$ weakly or strongly characteristic for $g$ if $g_k\cap
V^*\cdot S^{k-1}T^*\ot N\ne0$, resp. $g_k\cap S^kV^*\ot N\ne0$,
for $k=r_\text{max}(g)$.
 \end{dfn}
Note that now the notion of strong characteristicity is not
opposite to weak non-characteristicity unless the system has a
pure order $r_\text{min}=r_\text{max}$ (and the same for weak
characteristicity and strong non-characteristicity).

Relation between characteristic covectors and characteristic
subspaces were clarified by Guillemin \cite{G} in the case of
first order involutive systems. We extend his result to the case
of arbitrary orders in Theorem \ref{thm2}.

\section{\hps Involutivity}\label{S3}

 \abz
The classical Cartan's definition of involutivity involves
quasi-regular sequences. Namely a subspace $g_k\subset S^kT^*\ot
N$ is involutive if for some and hence for any generic basis
$v_1,\dots,v_n$ of $T$ the maps
 $$
\d_{v_i}: g_k^{(1)}\cap S^{k+1}\op{ann}\langle v_1,\dots,
v_{i-1}\rangle\ot N\to g_k\cap S^k\op{ann}\langle v_1,\dots,
v_{i-1}\rangle\ot N \eqno(*)
 $$
are surjective for all $1\le i\le n$.

It is well-known (see Serre's letter in \cite{GS}, also
\cite{BCG$^3$}) that this is equivalent to the requirement
 $$
H^{i,j}(g)=0 \text{ for }i\ge k, \eqno(**)
 $$
where $g$ is the system generated by $g_k$, and so ($g_k$
involutive)$\,\Rightarrow\,$($g_k^{(1)}$ involutive). Basing on
this homological characterization a system $g$ of pure order $k$
is called involutive if its generating subspace $g_k$ is such.

But there are several ways of generalizing this for arbitrary
symbolic systems:

 \begin{enumerate}
 \item[$I_1$:]
$H^{i,1}(g)=0\Longrightarrow H^{i,j}(g)=0\ \forall j>1$.
 \item[$I_2$:]
For $k\notin\op{ord}(g)$ and a generic basis $\{v_i\}$ the
following maps are surjective: $\d_{v_i}: g_k\cap
S^k\op{ann}\langle v_1,\dots, v_{i-1}\rangle\ot N\to g_{k-1}\cap
S^{k-1}\op{ann}\langle v_1,\dots, v_{i-1}\rangle\ot N$.
 \item[$I_3$:]
Denoting $\op{ord}(g)=\{r_1<\dots<r_s\}$, there is a splitting
$T=\oplus_{j=1}^sU_j$ and a basis $\{v_i\}$ subordinated to it
such that the above maps $\d_{v_i}$ are surjective unless $k=r_m$
is an order and $i$ corresponds to $v_i\in U_m$.
 \end{enumerate}

One easily proves $I_3\Rightarrow I_2\Rightarrow I_1$, but the
implications are irreversible. Properties $I_1,I_2$ are too weak
for the general definition of involutivity and $I_3$ seems to be
too strong (this property holds for the direct sum $g\subset
S(\sum T_i^*)\ot(\sum N_i)$ of involutive systems $g_i\subset
ST_i^*\ot N_i$). So we give:

 \begin{dfn}\po
A symbolic system $g\subset ST^*\ot N$ is called involutive if
each subspace $g_k\subset S^kT^*\ot N$ is involutive. When
$\op{ord}(g)=\{r_1,\dots,r_s\}$, this is a condition only for
$k=r_i$, $1\le i\le s$.
 \end{dfn}

This definition most appropriately reflects the dual picture of
quasi-regular sequences in the symbolic module $g^*$ (\cite{GS,
BCG$^3$}) known for pure order systems (it is also interesting to
investigate involutivity coupled with property $I_3$).

Let us denote by $g^{|k\rangle}$ the symbolic system generated by
all differential corollaries of the system deduced from the order
$k$:
 $$
g^{|k\rangle}_i=\left\{\begin{array}{ll}S^iT^*\ot N,&\text{ for }i<k;\\
g_k^{(i-k)},&\text{ for }i\ge k.\end{array}\right.
 $$

 \begin{theoreme}\po
A system $g$ is involutive iff $H^{i,j}(g^{|k\rangle})=0$ for all
$i\ge k$ (this condition is to be checked for $k\in\op{ord}(g)$
only).
 \end{theoreme}

 \begin{proof}
This follows from the classical equivalence
$(*)\,\Leftrightarrow\,(**)$ because involutivity of $g$ means
involutivity of the pure order systems
$g^{|r_1\rangle},\dots,g^{|r_s\rangle}$.%
 \end{proof}

In particular, $H^{i,j}(g)=0$ for $i\notin\op{ord}(g)-1$,
$(i,j)\ne(0,0)$, and properties $I_1,I_2$ follow from
involutivity. This however is not invertible:

 \begin{examp}\po\rm
Consider the system $u_{xx}=0,u_{yy}=0,v_{yyy}=0$. The only
non-zero cohomologies are $H^{0,0}(g)=\R^1$, $H^{1,1}(g)=\R^2$,
$H^{2,1}(g)=\R^1$, $H^{2,2}(g)=\R^1$. So $I_1$ holds. The
surjectivity requirement $I_2$ holds as well. But the pure order 2
system $g^{|2\rangle}$ is not involutive because
$H^{2,2}(g^{|2\rangle})=\R^1$. Thus $g$ is not involutive.
 \end{examp}

 \begin{examp}\po\rm
The system $g$ given by equations $u_{xx}=0,u_{xy}=0,u_{yyz}=0$ is
involutive. But it does not satisfy the property $I_3$.
 \end{examp}

\section{\hps Proof of Theorem \ref{thm1}}\label{S4}

 \abz
Our proof is based on a spectral sequence constructed in
\cite{KL$_2$} for the need of a reduction theorem (which means
that instead of projection of the symbolic system $g$ to $SW^*\ot
N$ we intersect it with $SV^*\ot N$).

Define a filtration in the $l$-th Spencer complex, induced by the
filtration in $\La T^*$ via the powers of $V^*$:
 $$
F^{p,q}=g_{l-p-q}\ot\La^pV^*\we\La^qT^*.
 $$

 \begin{lem}\po
The filtration is monotone decreasing, $F^{p+1,q-1}\subset
F^{p,q}$, and is preserved by the $\d$-map, $\d F^{p,q}\subset
F^{p,q+1}$. \qed
 \end{lem}

This filtration determines the spectral sequence of Leray-Serre type with
 $$
E_0^{p,q}=F^{p,q}/F^{p+1,q-1}=g_{l-p-q}\ot\La^pV^*\ot\La^qW^*.
 $$
The differential $d_0:E_0^{p,q}\to E_0^{p,q+1}$ acts by $W$ and so
 $$
E_1^{p,q}=H^{l-p-q,q}(g,\d')\ot\La^pV^*,
 $$
where $\d'$ is the induced differential (along $W$).

Denote $\Upsilon^{i,j}=V^*\cdot\d'(S^iT^*\ot N\ot\La^{j-1}W^*)$,
 $\Xi^{i,j}=\d'(g_{i+1}\ot\La^{j-1}W^*)\cap \Upsilon^{i,j}$.

If $V^*$ is strongly non-characteristic, then $\Xi^{i,1}=0$. In
fact, suppose that for some $p\in g_{i+1}$ we have: $\d'p\in
\Upsilon^{i,1}$. Then $(\d'p)|_W=0$ and so $p\,|_W=0$. This
contradicts injectivity of the projection $g_{i+1}\to\tg_{i+1}$.
However for $j>1$ the term $\Xi^{i,j}$ can be non-zero.

 \begin{lem}\po\label{lm5}
Let $g$ be a symbolic system, $V^*$ be strongly non-characteristic
and $k=r_\text{min}(g)$. The cohomology of $g$ with respect to the
induced differential $\d'$ are related to the Spencer cohomology
of the restricted symbolic system as follows:
 $$
H^{i,j}(g,\d')=\left\{
\begin{array}{ll}
0,&\!\! i<k-1,j>0;\\
S^iV^*\ot N,&\!\! i\le k-1,j=0;\\
H^{i,j}(\tg,\bar\d)\oplus \Upsilon^{i,j}/{\Xi^{i,j}},&\!\! i=k-1,j>0;\\
H^{i,j}(\tg,\bar\d)/\Xi^{i-1,j+1},&\!\! i=k;\\
H^{i,j}(\tg,\bar\d),&\!\! i>k.
\end{array} \right.
 $$
 \end{lem}

The third and forth lines above represent the cohomology
non-canonically. We actually mean here the exact sequences:
 \begin{gather*}
0\to\Upsilon^{k-1,j}/{\Xi^{k-1,j}}\to H^{k-1,j}(g,\d')\to
H^{k-1,j}(\tg,\bar\delta)\to0,
\\
0\to H^{k,j-1}(g,\d')\to H^{k,j-1}(\tg,\bar\delta)\to\Xi^{k-1,j}
\to0.
 \end{gather*}

 \begin{proof}
The restriction map induces an isomorphism of the complexes
 $$
  \begin{array}{ccccccccc}
0 & \to & g_l & \stackrel{\d'}\longrightarrow & g_{l-1}\ot W^* &
\stackrel{\d'}\longrightarrow & g_{l-2}\ot\La^2W^*
& \stackrel{\d'}\longrightarrow & \dots\\
&& \downarrow && \downarrow && \downarrow\\
0 & \to & \tg_l & \stackrel{\bar\d}\longrightarrow & \tg_{l-1}\ot
W^* & \stackrel{\bar\d}\longrightarrow & \tg_{l-2}\ot\La^2W^* &
\stackrel{\bar\d}\longrightarrow & \dots
  \end{array}
 $$
at first $l-k+1$ terms and hence an isomorphism of cohomologies.
For $i<k$ we have: $g_i=S^iT^*$, $\tg_i=S^iW^*$. So the boundary
cohomologies $H^{k,j-1}$, $H^{k-1,j}$ make the only difference,
being found from the commutative diagram:
 $$
  \begin{array}{ccccccc}
\dots & \to & g_k\ot\La^{j-1}W^* & \stackrel{\d'}\longrightarrow &
S^{k-1}T^*\ot N\ot\La^j W^* & \to & \dots\\
&& \wr\!\downarrow && \downarrow \\
\dots & \to & \tg_k\ot\La^{j-1}W^* & \stackrel{\bar\d}
\longrightarrow & S^{k-1}W^*\ot N\ot\La^j W^* & \to & \dots
  \end{array}
 $$
\vspace{-40pt}

 \end{proof}

 \begin{lem}\po\label{lem6}
If $g$ is involutive, then the differentials
$d_r^{p,q}:E_r^{p,q}\to E_r^{p+r,q-r+1}$ are trivial for $r>0$,
save for the map $d_1^{p,0}$ with $l-p<k=r_\text{min}(g)$, which
is the $\d$-differentiation along $V$.
 \end{lem}

 \begin{proof}
We prove at first this statement for the case of pure order $k$
symbolic system $g$. Afterwards we deduce the general case.

Let $g$ be involutive. We prove by induction on $l\ge k$ that
$H^{i,l-i}(g,\d')$ vanishes except for $i=k-1$. Let $\dim V=t$.

Then the table $E_1^{p,q}$ consists of $(t+1)$ columns:
$E_1^{0,q}=H^{l-q,q}(g,\d')$, $E_1^{1,q}=H^{l-1-q,q}(g,\d')\ot
V^*$, $\dots$, $E_1^{t,q}=H^{l-t-q,q}(g,\d')\ot \La^tV^*$. By the
induction hypothesis for $0<i\le t$ there is only one non-zero
term among $E_1^{i,q}$, $q>0$, corresponding to $q=l-i-k+1$ (if
this number is non-negative, otherwise all terms vanish).

Also the row $(E_1^{p,0},d_1^{p,0})$ is exact except for the left
boundary position, whence $E_2^{p,0}=0$ for $p\ne l-k+1$. For
$p=l-k+1\le t$ we have:
$E^{p,0}_2=H(E_1^{p,0},d_1^{p,0})=\d(S^kV^*\ot N\ot\La^{l-k}V^*)$.

Thus since only the elements of anti-diagonal $p+q=l-k+1$ can
survive in $E_\infty$, the single non-zero term in the $E_1^{0,q}$
column except $E_1^{0,l-k+1}=H^{k-1,l-k+1}(g,\d')$ can be
$E_1^{0,l-k}=H^{k,l-k}(g,\d')$ provided that one of the
differentials $d_i^{0,l-k}$, $i=1,\dots,t$, is injective.

But then in the spectral sequence for $(l+t)$-th Spencer complex
we find a non-zero term $E_1^{t,l-k}=H^{k,l-k}(g,\d')\ot \La^tV^*$
and to kill its contribution to the Spencer group $H^{k,*}(g)=0$
(involutivity) we need to assume a non-zero term
$E_1^{i,l+t-k-i-1}=H^{k+1,l+t-k-i-1}(g,\d')\ot\La^iV^*$, $0\le
i<t$. Continuing we obtain an infinite sequence of non-zero groups
$H^{s,q_s}(g,\d')=H^{s,q_s}(\tg)$, $s\to\infty$. But this
contradicts Poincar\'e $\d$-lemma, according to which $\dim
H^{*,*}(\tg)<\infty$.

Now consider the general case, when the symbolic system has
different orders. To prove vanishing of $d_r^{p,q}$ for $r>0$ we
use an alternative (but standard) definition of the $r$-th term of
the spectral sequence:
 $$
E_r^{p,q}=Z_r^{p,q}/(Z_{r-1}^{p+1,q-1}+B_{r-1}^{p,q}),
 $$
where $Z_r^{p,q}=\{\oo\in F^{p,q}\,|\,\d\oo\in F^{p+r,q-r+1}\}$
are $r$-th order cocycles and $B_r^{p,q}=\{\d\te\in
F^{p,q}\,|\,\te\in F^{p-r,q+r-1}\}$ are $r$-th order coboundaries.
Since the differential $d_r^{p,q}$ is induced by $\d$ it is
helpful to introduce the following spaces (we stress dependence on
the symbolic system for these terms):
 $$
\hat E_r^{p,q}(g)=Z_r^{p,q}/(Z_{r-1}^{p+1,q-1}+\op{Ker}\d\cap
Z_r^{p,q}),\qquad \check E_r^{p,q}(g)=Z_r^{p,q}/B_{r-1}^{p,q}.
 $$
Actually, the differential $d_r^{p,q}$ factorizes via the natural
$\d$-induced map $\bar d_r^{p,q}:\hat E_r^{p,q}(g)\to \check
E_r^{p+r,q-r+1}(g)$ and the natural projections
$\varsigma_r^{p,q}:\check E_r^{p,q}(g)\to E_r^{p,q}(g)$ and
$\varrho_r^{p,q}:E_r^{p,q}(g)\to\hat E_r^{p,q}(g)$ as follows:
 $$
d_r^{p,q}=\varsigma_r^{p+r,q-r+1}\circ\bar d_r^{p,q}\circ
\varrho_r^{p,q}.
 $$

Note that $Z_r^{p,q}(g)$ as well as $Z_{r-1}^{p+1,q-1}(g)$ depends
only on the graded term $g_{l-p-q}$, if we consider the $l$-th
Spencer complex. Thus $\hat E_r^{p,q}(g)=\hat
E_r^{p,q}(g^{|l-p-q\rangle})$ and consequently obtain the
following commutative diagram:
 $$
 \begin{array}{ccccccc}
\hspace{-3pt} E_r^{p,q}(g) & \to & E_r^{p,q}(g^{|l-p-q\rangle}) &
\stackrel{d_r^{p,q}}\longrightarrow &
E_r^{p+r,q-r+1}(g^{|l-p-q\rangle}) & \to &
E_r^{p+r,q-r+1}(g)\\
\hspace{-3pt} \downarrow{\scriptstyle \varrho_r} &&
\downarrow{\scriptstyle \varrho_r} && \uparrow{\scriptstyle
\varsigma_r} &&
\uparrow{\scriptstyle \varsigma_r}\\
\hspace{-3pt} \hat E_r^{p,q}(g) & \stackrel{\sim}\to & \hat
E_r^{p,q}(g^{|l-p-q\rangle}) & \stackrel{\bar
d_r^{p,q}}\longrightarrow & \check
E_r^{p+r,q-r+1}(g^{|l-p-q\rangle}) & \to & \check
E_r^{p+r,q-r+1}(g).
 \end{array}
 $$
In this diagram all arrows except differentials $d_r,\bar d_r$ are
projections and thus the differential $d_r^{p,q}(g)$ for the
general symbolic system $g$ factorizes via the differential
$d_r^{p,q}$ for the pure order involutive system
$g^{|l-p-q\rangle}$.

Since $d_r^{p,q}=0$, $q,r>0$, for pure order involutive systems as
was proved above, we obtain the same conclusion $d_r^{p,q}(g)=0$
(in the same range) for the general involutive symbolic system
$g$.
 \end{proof}

 \begin{lem}\po\label{lem6,5}
Let $\tg$ be an involutive system of pure order $k$. Then for
$r>0$ the differentials $d_r^{p,q}$ vanish as in Lemma \ref{lem6}
(except $d_1^{p,0}$, $l-p<k$) and $g$ is also involutive of pure
order $k$.
 \end{lem}

 \begin{proof}
If $\tg$ is involutive, then $E_1^{p,q}$ has support in the lines
$p+q=l-k+1$ and $q=0$. The latter is exact when equipped with the
differential $d_1$, except for the term $E_1^{p,0}$, $p=l-k+1\le
t$, and the former survives until $E_\infty$ by the graphical
evidence.

Now involutivity of $g$ follows from Lemma \ref{lm5} and the
spectral sequence, since $H^{i,j}(\tg)=0$ for $i\ne k-1$ implies
$H^{i,j}(g)=0$ for $i\ne k-1$.
 \end{proof}

 \begin{lem}\po\label{lm7}
Let the system $g$ be involutive. Then the system $\tg$ is also
involutive and \,$\Xi^{i,j}=0$ for $i\ge r_\text{min}(g)-1$.
 \end{lem}

 \begin{proof}
By our definition the first statement suffices to prove for pure
order $k$ systems. For $k=1$ the claim that $\tg$ is involutive is
the Guillemin's Theorem A and for $k>1$ it follows by the
equivalence reduction, see the end of \S\ref{S5}.

Since $H^{i,j-1}(\tg)=0$ for $i\ge k$, the homomorphic image of
this group $\Xi^{i-1,j}$ vanishes too. The same obviously holds
for a general system $g$ and $k=r_\text{min}(g)$, since the term
$\Xi^{i-1,j}$ depends only on the subspace $g_i$ and hence on the
involutive system $g^{|i\rangle}$.
 \end{proof}

Now we can finish the proof of Theorem \ref{thm1}. By Lemma
\ref{lem6} we have: $E_\infty^{p,q}=E_1^{p,q}$ for $q>0$,
$E_\infty^{p,0}=E_2^{p,0}=0$ for $p\ne l-k+1$,
$k=r_\text{min}(g)$, and $E_2^{l-k+1,0}=\Pi^{k-1,l-k+1}$. Thus we
conclude for $l>0$:
 $$
H^{l-j,j}(g,\d)\simeq\oplus_{p+q=j}\,E_\infty^{p,q}=
\oplus_{q>0}\, H^{l-j,q}(g,\d')\ot\La^{j-q}V^*\oplus
\d_{r_\text{min}(g)}^{l-j+1}\cdot\Pi^{l-j,j}
 $$
and the claim follows from Lemmata \ref{lm5} and \ref{lm7}.
\qed\vspace{5pt}

\section{\hps Proof of the Corollary}\label{S4.5}

 \abz
For $r_\text{min}(g)=1$ the sequence of the Corollary is exact for
all $i\ge0$. When $i=0$ it reads:
 \begin{multline*}
0\to S^jV^*\ot H^{0,0}(g)\to S^{j-1}V^*\ot H^{0,1}(g)\to\dots\\
\dots\to V^*\ot H^{0,j-1}(g) \to H^{0,j}(g)\to H^{0,j}(\tg)\to 0.
 \end{multline*}
This follows from the formula of Theorem \ref{thm1}, which in the
considered case can be rewritten as:
 $$
H^{0,j}(g)=\oplus_{q\ge0}H^{0,q}(\bar g)\ot\La^{j-q}V^*.
 $$
Substitution of this into the above sequence decomposes it into
the sum of the trivial Spencer complexes
$(S^{\a-t}V^*\ot\La^tV^*,\d)$ tensorially multiplied with
$H^{0,s}(\bar g)$.

Of course, the decomposition is not natural, so this argument is
not justified. But we can filter the cohomology $H^{0,j}(g)$ via
the spectral sequence:
 \begin{multline*}
H^{0,j}(g)=F_\infty^{0,j}\supset F_\infty^{1,j-1}\supset\dots
\supset F_\infty^{j,0}\ \text{ with }\\
F_\infty^{a,b}/F_\infty^{a+1,b-1}=E_\infty^{a,b}=H^{0,b}(\tg)\ot\La^aV^*.
 \end{multline*}
The associated graded sum is as in the considered formula.
Therefore we can filter the above complex and the consecutive
quotients are exact. The required exactness of the whole complex
follows.

 \begin{rk}\poo
In \cite{G} involutivity for pure order 1 restricted systems was
deduced from the exactness of the above sequence. Here we use the
opposite idea, concluding exactness from a by-product (or tool) of
Theorem \ref{thm1} on involutivity.
 \end{rk}

For $i=r_\text{min}(g)-1>0$ we have the following complex:
 \begin{multline*}
0\to S^{j,i}\ot N\to S^{j-1}V^*\ot H^{i,1}(g)\to\dots\\
\dots\to V^*\ot H^{i,j-1}(g) \to H^{i,j}(g)\to
H^{i,j}(\tg)\oplus\Upsilon^{i,j}\to 0.
 \end{multline*}
This complex is again filtered via the filtration of cohomology
 \begin{multline*}
H^{i,j}(g)=F_\infty^{0,j}\supset F_\infty^{1,j-1}\supset\dots
\supset F_\infty^{j,0}\ \text{ with }\\
F_\infty^{a,b}/F_\infty^{a+1,b-1}=E_\infty^{a,b}=(H^{i,b}(\tg)\oplus
\Upsilon^{i,b})\ot\La^aV^*\oplus\d_0^b\cdot\Pi^{i,a}
 \end{multline*}
The consecutive quotient complexes equal
$(S^{\a-t}V^*\ot\La^tV^*,\d)$ tensorially multiplied with
$H^{\b,s}(\bar g)\oplus\Upsilon^{\b,s}$ (again this summation is
not natural, so one should consider in stead the short exact
sequence as after Lemma \ref{lm5} and perform an additional
factorization) and so are exact. The last occurring complex is:
 $$
0\to S^{j,i}\to S^{j-1}V^*\ot\Pi^{i,1}\to\dots\to
V^*\ot\Pi^{i,j-1}\to\Pi^{i,j}\to0.
 $$
Its exactness follows from the following anti-commutative diagram
(or bi-complex: the sum of compositions of arrows along the
boundary of a square is zero), in which rows and columns are
exact, save for the one-term sequences:
 $$
 \begin{array}{ccccc}
\ddots&\uparrow&&\ddots\\
\dots\to& S^aV^*\ot S^{b-1}V^*\ot\La^{c+1}V^* &\rightarrow&
S^{a-1}V^*\ot S^{b-1}V^*\ot\La^{c+2}V^* &\to\dots\\
& \uparrow && \uparrow &\vspace{3pt}\\
 \dots\to& S^aV^*\ot S^bV^*\ot\La^cV^* &\rightarrow& S^{a-1}V^*\ot
S^bV^*\ot\La^{c+1}V^* &\to\dots\\
& \uparrow && \uparrow &\vspace{3pt}\\
&\ddots&\rightarrow& S^{a-1}V^*\ot S^{b+1}V^*\ot\La^cV^*&\to\dots
 \end{array}
 $$

In the remaining cases $i\ge r_\text{min}(g)$ and the complex from
the corollary equals:
 $$
0\to S^{j-1}V^*\ot H^{i,1}(g)\to\dots\to V^*\ot H^{i,j-1}(g) \to
H^{i,j}(g)\to H^{i,j}(\tg)\to 0.
 $$
Again we have a decreasing filtration $F_\infty^{t,j-t}$ of
$H^{i,j}$ with
 $$
F_\infty^{a,b}/F_\infty^{a+1,b-1}=E_\infty^{a,b}=
H^{i,b}(\tg)\ot\La^aV^*.
 $$
Thus the considered complex is filtered with all the consecutive
quotients being exact and the claim follows.

\section{\hps Proof of Theorem \ref{thm2}}\label{S5}

 \abz
We give an indirect proof, though a direct approach, similar to
the one presented in Appendix \ref{S9} and using the Corollary, is
plausible.

Recall that every system of PDEs of higher orders can be
equivalently written as a system of first order equations. This is
achieved via the map $\mathcal{E}\subset J^k(\pi)\hookrightarrow
J^1(J^{k-1}(\pi))$. Let us call this composition map the
equivalence reduction (er).

 \begin{examp}\po\ \rm\label{ex_char}
The equation $u_{xy}=0$ on the plane $T=\R^2(x,y)$ is equivalent
to the following system of the first order: $p_y=0$, $q_x=0$
($p=u_x,q=u_y$). We will identify $T$ with its tangent spaces and
consider the corresponding symbolic system $g$ on $T$. Let $W$ be
a proper subspace of $T$ which equals neither $\R^1(x)$ nor
$\R^1(y)$. The corresponding subspace $V^*=\op{ann}(W)\subset T^*$
is weakly characteristic for the second order PDE, but is not for
the equivalent first order system.

Thus the notion of weakly characteristic subspace is not invariant
under the equivalence reductions. However we will show the strong
characteristicity is well-posed.

Similarly, if $g$ is a scalar ($\dim N=1$) symbolic system on $T$
with $\dim T>1$ generated by one higher order PDE, then any
1-dimensional subspace $V^*$ is weakly characteristic, while it is
weakly (in this case also strongly) characteristic for the first
order reduction of $g$ iff the corresponding covector is
characteristic.
 \end{examp}

On the algebraic level the above equivalence reduction is obtained
via the embedding $\d:S^kT^*\to S^{k-1}T^*\ot T^*$, which induces
the following correspondence:
 $$
S^kT^*\ot N\supset g_k\ \rightsquigarrow\
\hat{g}_1=\op{er}_k(g_k)\subset T^*\ot(S^{k-1}T^*\ot N).
 $$

More generally for $l\ge k$ the coupling map $S^{k-1}T\ot
S^lT^*\to S^{l-k+1}T^*$ yields:
 $$
S^lT^*\ot N\supset g_l\ \rightsquigarrow\
\hat{g}_{l-k+1}=\op{er}_k(g_l)\subset
S^{l-k+1}T^*\ot(S^{k-1}T^*\ot N).
 $$
The map $\op{er}_k$ acts on elements as follows:
 $$
\prod_{i=1}^m a_i\ot\x\mapsto\sum_{i_1,\dots,i_{k-1}}
\dfrac{l!}{(l-k+1)!}\prod_{j\ne i_s} a_j\ot(a_{i_1}\cdots
a_{i_{k-1}})\ot\xi
 $$

We shall show this correspondence is respected by the prolongation
procedure, so that it descends to symbolic systems with
$r_\text{min}(g)=k$.

 \begin{lem}\po
The subspaces $\op{er}_k(g_l^{(i)}),\op{er}_k(g_l)^{(i)}\subset
S^{l+i-k+1}T^*\ot(S^{k-1}T^*\ot N)$ coincide for $l\ge k$.
 \end{lem}

 \begin{proof}
The following diagram commutes:
 $$
  \begin{array}{ccc}
S^lT^*\ot N\ot\La^{j-1}T^* &\rightsquigarrow &
S^{l-k+1}T^*\ot(S^{k-1}T^*\ot N)\ot\La^{j-1}T^*\\
\d\downarrow & & \d\downarrow \\
S^{l-1}T^*\ot N\ot\La^jT^* &\rightsquigarrow &
S^{l-k}T^*\ot(S^{k-1}T^*\ot N)\ot\La^jT^*.
  \end{array}
 $$
This suffices to check on generators:
 $$
  \begin{array}{ccc}
z_1^l\ot\xi\ot z_2\we\dots\we z_j &\stackrel{\op{\,er_k}}\mapsto &
\frac{l!}{(l-k+1)!}z_1^{l-k+1}\ot z_1^{k-1}\ot\xi\ot z_2\we\dots\we z_j\\
\downarrow & & \downarrow \\
lz_1^{l-1}\ot\xi\ot z_1\we z_2\we\dots\we z_j
&\stackrel{\op{\,er_k}}\mapsto & \frac{l!}{(l-k)!}z_1^{l-k}\ot
z_1^{k-1}\ot\xi\ot z_1\we z_2\we\dots\we z_j.
  \end{array}
 $$
Thus the equivalence reduction and the prolongation commute.
 \end{proof}

Consequently we can define equivalence reduction of a symbolic
system $g$. Let's denote the reduction $\op{er}_k(g)$ of a
symbolic system $g$ by $\hat{g}$.

 \begin{prop}\po\label{invol}
Let $g$ be a symbolic system of minimal order $r_\text{min}(g)\ge
k$. Then involutivity of the system $g$ is equivalent to
involutivity of the system $\hat{g}$.
 \end{prop}

 \begin{proof}
Consider at first a symbolic system $g$ of pure order $\ge k$. We
have the following isomorphism of complexes for $l\ge k$:
 $$
  \begin{array}{ccc}
\vdots & & \vdots\\
\d\downarrow & & \d\downarrow \\
g_{l+1}\ot\La^{j-1}T^* &\rightsquigarrow &
\hat{g}_{l-k+2}\ot\La^{j-1}T^*\\
\d\downarrow & & \d\downarrow \\
g_l\ot\La^jT^* &\rightsquigarrow & \hat{g}_{l-k+1}\ot\La^jT^*\\
  \end{array}
 $$
Thus $H^{i,j}(g)=H^{i-k+1,j}(\hat g)$, if $i>k-1$. But we can also
consider the derived systems $g^{|s\rangle}$ instead of $g$ and
get the same conclusion. The claim follows.
 \end{proof}

 \begin{lem}\po\label{lem12}
$V^*\subset T^*$ is strongly characteristic for $g$ iff
it is such for $\hat g$.
 \end{lem}

 \begin{proof}
If $g_k\cap S^kV^*\ot N\ne0$, then clearly ${\hat g}_1\cap
V^*\ot(S^{k-1}V^*\ot N)\ne0$. On the other hand, if ${\hat
g}_1\cap V^*\ot(S^{k-1}T^*\ot N)\ne0$, then the pairing of
$S^{k-1}T\ot N^*$ and some element $p\in g_k$ takes values in
$V^*$, whence $g_k\cap S^kV^*\ot N\ne0$.
 \end{proof}

 \begin{prop}\po
The characteristic varieties of the systems $g$ and $\hat g$
coincide.
 \end{prop}

 \begin{proof}
Actually, the characteristic variety $\op{Char}^\C(g)$ is defined
by the characteristic ideal $I(g)=\op{ann}(g^*)\subset ST$
(\cite{S}). But from the description of this ideal given in
\cite{KL$_2$} we see that $I(g)=I(\hat g)$.
 \end{proof}

Thus we have reduced Theorem \ref{thm2} for pure order $k$ systems
to its partial case for $k=1$, i.e.\ Theorem B of \cite{G} (see
Appendix \ref{S9} for the proof). If $g$ is a general involutive
system, then $g^{|k\rangle}$, $k\le r_\text{max}(g)$, is a pure
order involutive system and already proved part of the statement
implies the whole claim. \qed

\vspace{4pt}
 \begin{Proof}{Proof that ($g$ involutive)$\,\Rightarrow$\,($\tg$ involutive)}
As noted in \S\ref{S4} this implication suffices to prove for pure
order $k$ systems. We do it by reducing to the case of first
order.

By Proposition \ref{invol} the equivalence reduction $\hat g$ is
an involutive first order system. Since $V^*$ is strongly
non-characteristic for $g$ it is also such for $\hat g$ (to this
side the claim is true: otherwise the subspace $V^*$ is strongly
characteristic for $\hat g$ and we apply Lemma \ref{lem12}). Thus
by Guillemin's Theorem A the reduction $\bar{\hat g}\subset
SW^*\ot(S^{k-1}T^*\ot N)$ is involutive.

Since $V^*$ is strongly non-characteristic, the coefficients
reduction $S^{k-1}T^*\ot N\to S^{k-1}W^*\ot N$ maps the system
$\bar{\hat g}$ isomorphically onto its image, which is the
equivalence reduction of the restricted system $\hat{\tg}\subset
SW^*\ot(S^{k-1}W^*\ot N)$. But the map $\op{eq}_k$ does not change
involutivity, and therefore the system $\tg$ is involutive.
 \end{Proof}

\section{\hps Examples}\label{S6}

 \abz
Here we demonstrate that all assumptions in our results are
essential. At first we consider the statement of Theorem
\ref{thm2}.

 \begin{examp}\po\ \rm\label{ex4}
Consider a system generated by a subspace $g_k\subset S^kT^*$ of
codimension 1. One easily checks the system $g$ is involutive
(this follows immediately from the reduction theorem of
\cite{KL$_2$}). The set of characteristic covectors forms a
hypersurface in $\mathbb{P}^\C T^*$ of degree $k$ and so a generic
covector is non-characteristic. However all 1-dimensional
subspaces $V^*$ are weakly characteristic if $k>1$ and $\dim T>2$
(while weakly non-characteristic subspaces are plentiful too).
 \end{examp}

 \begin{examp}\po\ \rm
For the Laplace equation on $T$ with $\dim T\ge2$ all $V^*$ of
dimension 2 are strongly characteristic. But there are no real
characteristics. Thus working over $\C$ is important.
 \end{examp}

Let $g\subset ST^*$ be a scalar system of {\em complete
intersection type} (\cite{KL$_2$}), which means that if
$\op{codim}(g)=t$ and $\dim T=n$, we have: $t\le n$ and
$\op{codim}\op{Char}^\C(g)=t$.

 \begin{examp}\po\ \rm\label{examp6}
Let $\dim T=n>2$ and a symbolic system $g\subset ST^*$ of order 2
be given by $n$ equations of complete intersection type. There are
strongly characteristic subspaces $V^*$ of dimension $>1$, but the
system is of finite type and hence is free of characteristics.
Note that for this system
 $$
H^{0,0}(g)\simeq\R^1,\ H^{1,1}(g)\simeq\R^n,\dots,
H^{i,i}(g)\simeq\R^{\binom ni},\dots, H^{n,n}(g)\simeq\R^1.
\eqno(\ddag)
 $$
and so it is not involutive.
 \end{examp}

Now we study the situation of Theorem \ref{thm1}. Consider a
general symbolic complete intersection system $g$.

 \begin{prop}\po\label{1910}
If a subspace $V^*\subset T^*$ is strongly non-characteristic for
$g_{r_\text{max}(g)}$, then $\dim V^*\le m(1)+1$. The equality can
be achieved only if the system has two orders
$\op{ord}(g)=\{1,k\}$ with multiplicities $m(1),m(k)$, satisfying
$n=m(1)+m(k)$ (finite type) and either $k=2$ or $n-m(1)=m(k)=1$.
 \end{prop}

Notice that if $m(1)\ne0$, then strongly non-characteristic
subspace $V^*$ for $g$ should satisfy the inequality $\dim V^*\le
m(1)$ and every such generic subspace is strongly
non-characteristic. However if $m(1)=0$ and $n>1$, then the
proposition implies that the system has pure order 2 and finite
type (provided it is a complete intersection).

 \begin{proof}
If $\op{ord}(g)=\{1\}$, the claim is obvious. So consider the case
with higher order equations $\op{ord}(g)=\{1<k_1<\dots<k_t=k\}$.
Note that equations of the first order can be normalized to be
$u_{x_{n-m(1)+1}}=0,\dots,u_{x_n}=0$ and these variables $x_i$ can
be excluded for considerations in higher orders.

Namely, we decompose $T=\tilde T\oplus U$, where $\tilde
T=\R^{n-m(1)}(x_1,\dots,x_{n-m(1)})$ and
$U=\R^{m(1)}(x_{n-m(1)+1},\dots,x_n)$ with $g_k\subset S^k\tilde
T^*$. Denote $\tilde V=V\cap\tilde T$. Then $g_k\cap\tilde
V^*\cdot S^{k-1}\tilde T^*=0$ and so $\sum
m(k_i)\ge\op{codim}(g_k\subset S^k\tilde T^*)\ge \dim \tilde
V^*\cdot S^{k-1}\tilde T^*$.

We have $\dim\tilde V\ge\dim V-m(1)$ and if $\dim\tilde V>0$, we
get:
 $$
n-m(1)\ge \sum_{i=1}^t m(k_i)\ge \op{codim}g_k\ge \dim\tilde
V^*\cdot S^{k-1}\tilde T^*\ge\dim \tilde V^*\cdot\tilde T^*\ge
n-m(1).
 $$
Thus we must have equalities everywhere and the claim follows.
 \end{proof}

Let us consider a scalar complete intersection system $g$, not of
the first order (everything is clear), where the equality in
Proposition \ref{1910} is achieved. Due to a remark before the
proof we then restrict to a pure second order finite type system
of complete intersection type.

We will use this example to show that Theorem \ref{thm1} does not
extend generally. The spectral sequence, having $E^{p,q}_1$ in a
product form until $p+q\ne l-k+1$ as in \S\ref{S4}, may seem to
split and stabilize, but we will show the differentials
$d_1^{p,q}\not\equiv0$.

 \begin{examp}\po\ \rm\label{examp7}
Let a scalar symbolic system $g$ be given by $g_0=\R^1$,
$g_1=T^*$, $g_2\subset S^2T^*$ of codimension $n$ and complete
intersection type, $g_{2+i}=g_2^{(i)}$. Then $\dim g_i=\binom ni$
and the non-zero cohomologies are listed in $(\ddag)$.

Let $V^*$ be a generic 1-dimensional subspace of $T^*$ and
$W=\op{ann}(V^*)$. Then
 $$
\tg_0=\R^1,\ \tg_1=W^*,\ \tg_2=S^2W^*,\ \tg_{2+i}\subset
S^{2+i}W^*\text{ has codimension }
\tbinom{n+i}{i+2}-\tbinom{n}{i+2}.
 $$

The first non-trivial case from our point of view is $n=3$ ($\dim
W=2$), which we consider in details (the case of arbitrary $n$ is
absolutely similar). In this case the only non-zero Spencer
$\d$-cohomologies are:
  \begin{gather*}
H^{0,0}(g,\d')=\R^1,\ H^{1,0}(g,\d')=V^*\simeq\R^1,\
H^{1,1}(g,\d')=V^*\ot W^*\simeq\R^2,\\
H^{2,1}(g,\d')=S^2V^*\ot W^*\simeq\R^2,\
H^{2,2}(g,\d')=S^2V^*\ot\La^2W^*\simeq\R^1,\\
H^{3,2}(g,\d')=S^3V^*\ot\La^2W^*\simeq\R^1.
 \end{gather*}
Thus the spectral sequences for $l$-th Spencer complex stabilize
by graphical reasons at $E_1$ for even $l$ ($=2,4,6$), converging
to $H^{i,l-i}(g)$.

For odd $l$ the following table describes the term $E_1^{p,q}$ (we
compactify the notations $S^iV^*=S^i_{V^*}$ etc):

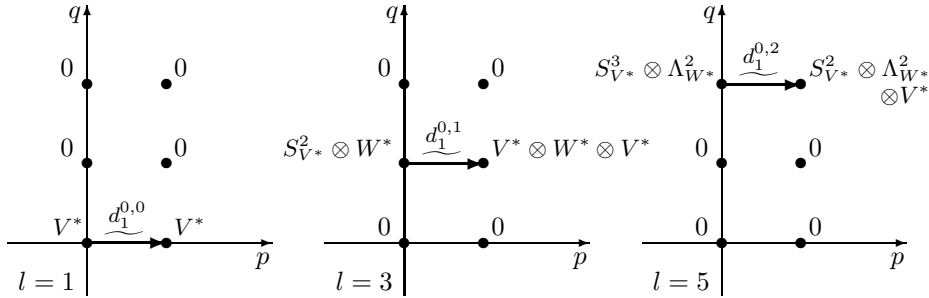
  \begin{picture}(440,130)
 \put(-10,30){\vector(1,0){100}}
 \put(20,10){\vector(0,1){110}}
 \put(20,30){\thicklines{$\vector(1,0){30}$}}
 \put(24,27){$\bf\widetilde{\hspace{20pt}}$}
 \put(28,38){\footnotesize $d_1^{0,0}$}
 \put(-6,12){$l=1$}
 \put(84,22){$p$}
 \put(13,115){$q$}
 \multiput(20,30)(0,30){3}{\circle*{4}}
 \multiput(50,30)(0,30){3}{\circle*{4}}
 \put(10,63){$0$}
 \put(53,63){$0$}
 \put(10,93){$0$}
 \put(53,93){$0$}
 \put(7,33){\small $V^*$}
 \put(53,33){\small $V^*$}
 \put(110,30){\vector(1,0){100}}
 \put(140,10){\vector(0,1){110}}
 \put(140,60){\thicklines{$\vector(1,0){30}$}}
 \put(144,57){$\bf\widetilde{\hspace{20pt}}$}
 \put(148,68){\footnotesize $d_1^{0,1}$}
 \put(114,12){$l=3$}
 \put(204,22){$p$}
 \put(133,115){$q$}
 \multiput(140,30)(0,30){3}{\circle*{4}}
 \multiput(170,30)(0,30){3}{\circle*{4}}
 \put(130,33){$0$}
 \put(173,33){$0$}
 \put(130,93){$0$}
 \put(173,93){$0$}
 \put(94,63){\small $S^2_{V^*}\ot W^*$}
 \put(173,63){\small $V^*\ot W^*\ot V^*$}
 \put(230,30){\vector(1,0){100}}
 \put(260,10){\vector(0,1){110}}
 \put(260,90){\thicklines{$\vector(1,0){30}$}}
 \put(264,87){$\bf\widetilde{\hspace{20pt}}$}
 \put(268,98){\footnotesize $d_1^{0,2}$}
 \put(234,12){$l=5$}
 \put(324,22){$p$}
 \put(253,115){$q$}
 \multiput(260,30)(0,30){3}{\circle*{4}}
 \multiput(290,30)(0,30){3}{\circle*{4}}
 \put(250,33){$0$}
 \put(293,33){$0$}
 \put(250,63){$0$}
 \put(293,63){$0$}
 \put(212,93){\small $S^3_{V^*}\ot\La^2_{W^*}$}
 \put(293,93){\small $S^2_{V^*}\ot\La^2_{W^*}$}
 \put(320,83){\small $\ot V^*$}
  \end{picture}
Thus $E_2^{p,q}=0$ ($p+q=l$) and we get stabilization at $E_2$.
 \end{examp}

This example shows that Lemmata \ref{lem6}, \ref{lm7} are wrong in
the non-involutive case (namely: $d_1^{0,s}\ne0$ for $l=2s+1$ and
$\Xi^{1,2}\ne0$).

 \begin{examp}\po\rm
Let $T=N=\R^2=\langle\p_x,\p_y\rangle$. Consider the symbolic
system: $g_0=N$, $g_1=\op{so}(2)=\langle\p_y\ot dx-\p_x\ot dy
\rangle$, $g_{1+i}=g_1^{(i)}=0$ for $i>0$. The only non-zero
cohomology are: $H^{0,0}(g)=\R^2$, $H^{0,1}(g)=\R^3$,
$H^{1,2}(g)=\R$. Thus $g$ does not satisfy $I_1$ and so is not
involutive.

Let $W=\langle\p_y\rangle$. Its annulator $V^*=\langle dx\rangle$
is a non-characteristic subspace for the system $g$ (in this case
weakly and strongly) . The restriction to $W$ is equal to:
$\tg_0=N$, $\tg_1=\langle\p_x\ot dy\rangle$, $\tg_{1+i}=0$. Thus
$\tg$ is not involutive in the pure order sense, but it is
involutive with our definition of involutivity for multiple-order.
The Spencer cohomology equal: $H^{0,0}(\tg)=\R^2$,
$H^{0,1}(\tg)=\R^1$, $H^{1,1}(\tg)=\R^1$. So $\tg$ satisfies the
properties $I_1,I_2$, but does not satisfy $I_3$ (though it
satisfies a modified $I_3'$, where we split not only the base $T$,
but also the fiber $N$).

Therefore, we see that the second part of Theorem \ref{thm1} does
not hold for involutive system $\tg$ of arbitrary orders. Also the
formula of Theorem \ref{thm1} does not hold for the systems $g$
and $\tg$ (for instance for $i=j=1$).
 \end{examp}

\section{\hps Other results and a discussion}\label{S7}

 \abz
We deduce one more result from the spectral sequence of
\S\ref{S4}. A symbolic system $g$ is called $m$-acyclic if
$H^{i,j}(g)=0$ for $i\notin\op{ord}(g)-1$ and $0\le j\le m$.
Involutivity corresponds to the case $m=\dim T$.

 \begin{theoreme}\po
Let $V^*$ be strongly non-characteristic for the symbolic system
$g$ of pure order $k$. Then $g$ is $m$-acyclic iff $\tg$ is
$m$-acyclic.
 \end{theoreme}

In a weak form this also generalizes to general symbolic systems.

 \begin{proof}
The reasoning for the direct implication is the same as in the
proof of Theorem \ref{thm1} (see Lemma \ref{lem6}): We prove by
induction on $l\ge k$ that $H^{l-j,j}(g,\d')=0$,
$j\le\min\{m,l-k\}$. The base of induction is obvious. Let us
study at first the case $\dim V^*=1$.

Consider the spectral sequence $E_r^{p,q}$ of the $l$-th Spencer
complex. By induction hypothesis $E_1^{i,j}=0$ for all $i\ge1,j\le
m,l-k-1$. Thus all the differentials $d_r^{0,j}$ for $j<m,l-k$
vanish and since $H^{l-j,j}(g)=0$ for $j\le m,l-k$, we conclude
that the only non-zero group among $E_1^{0,j}$ for $j\le m,l-k$
can occur when $j=j_0=\min\{m,l-k\}$ and only when the
differential $d_1^{0,j_0}$ is injective.

So suppose $E_1^{0,j_0}=H^{l-j_0,j_0}(g,\d')\ne0$. Then for the
spectral sequence of $(l+1)$-st Spencer complex the group
$E_1^{1,j_0}=H^{l-j_0,j_0}(g,\d')\ot V^*\ne0$ and in order to have
$H^{l-j_0,j_0}(g)=0$ the group $E_1^{1,j_0}$ should be killed by
$d_1^{0,j_0}:E_1^{0,j_0}\to E_1^{1,j_0}$. Thus
$E_1^{0,j_0}=H^{l+1-j_0,j_0}(g,\d')\ne0$. Continuing this process
we obtain a sequence $H^{s-j_0,j_0}(g,\d')\ne0$, $s\to\infty$,
which cannot happen by the $\d$-lemma.

Thus the claim is proved for $\dim V=1$. When $\dim V=t>1$ we can
find a complete flag $\{0\}\subset V_1^*\subset\dots\subset
V_t^*=V^*$ of strongly non-characteristic subspaces and apply the
previous arguments successively. Since again the corresponding
terms $\Xi$ vanish, we have $H^{l-j,j}(\tg)=H^{l-j,j}(g,\d')=0$
for $j\le m,l-k$.

The reverse statement follows directly from the spectral sequence
of \S\ref{S4}.
 \end{proof}

An alternative approach to the direct implication is via the
equivalence reduction and Proposition 2 of \cite{GK}, which is
equivalent to our statement for the pure first order systems (or
via the long diagram chase as in \cite{G}, p. 275).

In particular, 2-acyclicity is very important since obstructions
for formal integrability belong to the groups $H^{i,2}(g)$. After
some number of prolongations the system becomes $m$-acyclic, even
involutive. The above result states that the places, where this
stabilization happens, is the same for the systems $g$ and $\tg$.

 \begin{rk}\poo
It is possible to consider vanishing of $\d$-cohomology on the
other part of the spectrum: $H^{i,j}(g)=0$ for $n-j\le m$. This
$m$-coacyclicity is not so wide-spread as $m$-acyclicity, but is
closer in spirit to the notion of involutivity (because it implies
existence of a quasi-regular sequence of length $m$, see also
Appendix \ref{S8}). Then by similar methods one proves for $i\ge
k$:
 $$
H^{i,j}(g)=0\ \forall j\ge n-m \Leftrightarrow H^{i,j}(\tg)=0\
\forall j\ge n-m.
 $$
 \end{rk}

Involutive systems became one of the most important classes of
PDEs during the profound investigation of differential equations
compatibility and integrability problem at the beginning of the
last century (\cite{C,J,V}). However Cohen-Macaulay systems, being
very important in the commutative algebra, were introduced into
differential equations context quite recently in \cite{KL$_2$}.
Our frequent example (\S\ref{S6}) of complete intersections is a
partial case.

Recall that projection to the r.h.s.\ in the formula $(\dag)$
determines the restriction $\tg$ of a symbolic system $g$, while
intersection with the l.h.s.\ yields the reduction $\ti g$. Recall
also (\cite{KL$_2$}) that a subspace $V^*$ is called transversal
if its complexification is transversal to the characteristic
variety $\op{Char}^\C(g)$.

Now we wish to compare these two classes of symbolic systems. The
following table shows an apparent duality between them:

  \begin{center}
 {\renewcommand{\arraystretch}{1.5}%
 \begin{tabular}{|p{5.5cm}|p{5.5cm}|}
\hline \hfil Involutive systems \hfil & \hfil
Cohen-Macaulay systems \hfil \\
 \hline
\hline Restriction $\tg$ of an involutive $g$ to a strictly
non-characteristic subspace $W\subset T$ is involutive & Reduction
$\ti g$ of a Cohen-Macaulay $g$ to a transversal subspace
$V^*\subset T^*$ is Cohen-Macaulay \\
\hline Restriction induces an isomorphism of symbolic systems
$g\simeq\tg$, but not $\d$-cohomologies $H^{i,j}(g)\ne
H^{i,j}(\tg)$ & Reduction yields an isomorphism of
$\d$-cohomologies $H^{i,j}(g)\simeq H^{i,j}(\ti g)$, but not
symbolic systems $g\ne\ti g$ \\
\hline The Spencer cohomology $H^{*,*}(g)$ is a free $\La
V^*$-module, and change of coefficients $SW\subset ST$ induces an
isomorphism of the module $g^*$ & The symbolic module $g^*$ is a
free $SW$-module, and coefficients change $\La V^*\subset\La
T^*$ induces an isomorphism of the cohomology $H^{*,*}(g)$ \\
\hline Provides a canonical Koszul resolution of the symbolic
module $g^*$ & Provides an effective calculation of the Spencer
groups
$H^{i,j}(g)$ \\
 \hline
 \end{tabular}}
  \end{center}

The statements gathered here are combinations of results from
\cite{G}, \cite{KL$_2$} and the present papers.

\appendix

\section{\hps Descended symbolic systems}\label{S8}

 \abz
Here we describe a new property of involutive systems, not related
to the main topic of the paper, but still important in our
discussion.

Given a subspace $g_k\subset S^kT^*\ot N$ we can consider its {\it
descended\/} subspace $\p g_k=\langle\d_vp\,|\,v\in T,p\in
g_k\rangle\subset S^{k-1}T^*\ot N$.

If $g=\{g_k\}$ is a symbolic system, then $\p g_k\subset g_{k-1}$.

 \begin{prop}\po
Let $n=\dim T$. Then $H^{i,n}(g)=0$ iff $g_i=\p g_{i+1}$.
 \end{prop}

 \begin{proof}
This follows from exactness of the sequence:
 $$
g_{i+1}\ot\La^{n-1}T^*\longrightarrow\p g_{i+1}\ot \La^nV^*\to0.
 $$

\vspace{-30pt}
 \end{proof}

We define the {\it descended\/} system $\p g$ by the rule: $(\p
g)_k=\p g_{k+1}$.

 \begin{lem}\po
The descender $\p g$ is a symbolic system.
 \end{lem}

 \begin{proof}
Let $v\in T$ and $p\in\p g_{k+1}$. Then $p=\sum\d_{w_i}q_i$ for
some $q_i\in g_{k+1}$, $w_i\in T$. We have: $\d_v
p=\sum\d_v\d_{w_i}q_i=\sum\d_{w_i}(\d_vq_i)\in\p g_k$.
 \end{proof}

 \begin{prop}\po
$H^{i,n}(g)=0$ for $i\ge k$ implies $H^{i-1,n}(\p g)=0$ for $i\ge
k$.
 \end{prop}

 \begin{proof}
By a Bourbaki's lemma for the symbolic module $g^*$ (see Serre's
letter in \cite{GS}) the first property is equivalent to the
existence of $v\in T$ such that $\d_v: g_{i+1}\to g_i$ is
epimorphic for $i\ge k$. But then $\d_v:\p g_{i+1}\to\p g_i$ is
epimorphic too and the claim follows.
 \end{proof}

Thus the descender $\p g$ characterizes vanishing of certain
cohomologies. Recall that this vanishing is closely related to the
existence of quasi-regular elements (\cite{GS}) and so to
involutivity. In fact, the descended system behave nicely w.r.t.
involutivity property:

 \begin{theoreme}\po
If $g$ is an involutive system, so is it descender $\p g$.
 \end{theoreme}

 \begin{proof}
The statement suffices to check for the systems of pure order $k$.
Consider the commutative diagram of $\d$-sequences:
 $$
 \begin{array}{ccccc}
\ddots&\vdots&&\vdots\\
\dots\to& g_{k+2}\ot\La^iT^*\ot\La^{n-2}T^* &\longrightarrow&
g_{k+1}\ot\La^{i+1}T^*\ot\La^{n-2}T^* &\to\dots\\
& \downarrow && \downarrow &\vspace{3pt}\\
\dots\to& g_{k+1}\ot\La^iT^*\ot\La^{n-1}T^* &\longrightarrow&
g_k\ot\La^{i+1}T^*\ot\La^{n-1}T^* &\to\dots\\
& \downarrow && \downarrow &\vspace{3pt}\\
\dots\to& \p g_{k+1}\ot\La^iT^*\ot\La^nT^* &\longrightarrow&
\p g_k\ot\La^{i+1}T^*\ot\La^nT^* &\to\dots\\
& \downarrow && \downarrow &\vspace{3pt}\\
& 0 & & 0 &
 \end{array}
 $$
The columns are exact (in the symbolic grading $\ge k$). The rows
are also exact, save for the bottom one. So this latter is exact
too (in the symbolic grading $\ge k-1$) by the standard diagram
chase.
 \end{proof}

Though $\p g\subset g$, the equality does not hold even for
involutive systems. For instance, the Frobenius type system
$g_{k-1}=S^{k-1}T^*\ot N$, $g_k=0$, exhibits the strict inclusion.
However, it is the only finite type involutive system of pure
order and one can expect that modulo sums with such systems other
involutive systems enjoy the considered property.

In any case, either $\p g=g$ or we have a canonical sub-system,
which can be considered as an intermediate integral of the system.
Such integrals are known to facilitate the integration
(\cite{KL$_2$}).

We can iterate the procedure and construct further descended
systems: $g\supset \p g\supset\p^2 g\supset\dots$. This systems of
strict inclusions necessary terminates and in a finite number of
steps we get the least descender $\p^\infty g$.

\section{\hps Guillemin's theorem on characteristics}\label{S9}

 \abz
For completeness we provide here a proof of Theorem B from
\cite{G}%
, following the author's ideas, but with a different exposition
and the coordinate-free approach. Namely we prove that a
characteristic subspace contains a characteristic covector (over
$\C$), the inverse statement being obvious.

Let $g$ be a first order symbolic system and $V^*$ a
characteristic subspace, which in this case is the same as a
strongly characteristic subspace. Let $V_0^*$ be its maximally
non-characteristic subspace. Shrinking $V^*$ if necessary we can
assume that codimension of $V_0^*$ in $V^*$ is one.

Let $\oo\in V^*\setminus V_0^*$. We assume also that $\oo$ is not
a characteristic covector, for otherwise we are done. Denote
$N'=[\oo\ot N\cap(V_0^*\ot N+g)]/\oo$. Since $V_0^*\ot N\cap g=0$,
we get a well-defined map:
 $$
\l:N'\to V_0^*\ot N,\qquad \x\mapsto\oo\ot\x\mod g.
 $$

 \begin{prop}\po
{\rm (i)} $\op{Im}(\l)\subset V_0^*\ot N'$.\\
{\rm (ii)} The following $\l$-generated sequence is a complex:
 $$
0\to N'\stackrel{\l}\to V_0^*\ot N'\stackrel{\l}\to \La^2V_0^*\ot
N'\to\dots
 $$
 \end{prop}

 \begin{proof}
(i) Let $\l_v(\x)=\langle\l(\x),v\rangle$, $v\in
V_0=\op{ann}(\oo)\subset V$, be the map $\l_v:N'\to N$. We have:
 $$
 \begin{array}{rcl}
& V_0^*\ot T^*\ot N &\\
& \rho\downarrow &\\
0\to S^2V_0^*\ot N\stackrel{\vp_0}\longrightarrow & V_0^*\ot
H^{0,1}(g)&\stackrel{\vp_1}\longrightarrow H^{0,2}(g)
 \end{array}
 $$
The vertical map is epimorphic and for $\oo\ot\l(\x)\in V_0^*\ot
T^*\ot N$ we have: $\vp_1\circ\rho(\oo\ot\l(\x))=0$. From
exactness of the horizontal complex (Quillen's Theorem, see our
Corollary) we get: $\oo\ot\l_v(\x)\in i_v S^2V_0^*\ot
N=V^*_0\ot N\mod g$, i.e.
$\l_v(\x)\in N'$ $\forall v\in V_0$.\\
(ii) This claim follows from the identity $\oo\we\oo=0$.
 \end{proof}

 \begin{lem}\po
There exists a non-zero vector $\x_0\in N'$ such that $\l(\x_0)\in
V_0^*\ot\x_0$.
 \end{lem}

 \begin{proof}
This equivalently means that the set of commuting  (by (ii) above)
linear operators has a common eigenvector (over $\C$).
 \end{proof}

Now $\oo\ot\x_0=p\ot\x_0\mod g$ for some covector $p\in V_0^*$,
i.e.\ the covector $\bar p=\oo-p\in V^*$ is characteristic: $\bar
p\in\op{Char}^\C(g)$.


\vspace{-12pt}
\hspace{-20pt} {\hbox to 12cm{ \hrulefill }}
\vspace{-1pt}

{\footnotesize
\hspace{-10pt}
Institute of Mathematics and Statistics, University of Troms\o, Troms\o\
90-37, Norway.

\hspace{-10pt}
E-mails: \quad kruglikov\verb"@"math.uit.no,
\quad lychagin\verb"@"math.uit.no.}
\vspace{-1pt}

\end{document}